
Content-Length: 15319

\documentstyle{amsppt}
\magnification=1200
\baselineskip=18pt
\NoBlackBoxes
\def\P{\Bbb P}

\def\ls{\vskip.25in}

\centerline{\bf On the quantum cohomology of the plane, old and new}
\ls
\centerline{Z. Ran}
\ls
\centerline{\bf Abstract}
{\it We describe a method for counting maps of curves of given genus (and
variable
moduli) to $\Bbb P^2$, essentially by splitting the $\Bbb P^2$ in two; then
specialising to the case of genus 0 we show that the method of quantum
cohomology may be viewed as the 'mirror' of the former method where one splits
the $\Bbb P^1$ rather than the $\Bbb P^2$, and we indicate a proof of the
associativity of quantum multiplication based on this idea.}
\ls
Recent work on Mirror Symmetry and Quantum Cohomology has contributed to a
revival
of interest in problems of a classical nature in Enumerative Geometry (cf. [F]
and
references therein).  These problems involve (holomorphic) maps
$$
f: C \to X
\tag1
$$
where $X$ is a fixed variety and $C$ is a compact Riemann surface whose moduli
are
sometimes fixed (`Gromov-Witten') but here will not be,
unless otherwise stated.  While the case $\dim X = 1$ is
not entirely without interest (cf. [D]), the problem begins in earnest with
$\dim X =2$ and
naturally the simplest such $X$ is $\Bbb P^2$.  Here the problem specifically
is to
count the images $f(C)$ of maps (1) where $C$ has genus $g,f(C)$ has degree $d$
and passes through
$3d + g-1$ fixed points in $\Bbb P^2$.  This problem has already, in essence,
been solved
in the author's earlier paper [R] by means of a recursive method (we note
however that the formula in [R],
(3c.1), (3c.3) is trivially misstated and the factor $c(\tilde{K}_1,
\tilde{K}_2)$
should not be present).

Our purpose here is twofold. In Sect.1 we
give a partial exposition of the method of [R] and illustrate it on
a couple of new examples, namely the curves of degree $d$ and genus
$g=\frac{(d-1)(d-2)}{2} - 2$
(i.e. with 2 nodes); and the rational quartics.  We recover classical formulae
due,
respectively, to Roberts [Ro] and Zeuthen [Z].  Hopefully, this will help make
the method
of [R] more accessible. In Sect.2 we show that the method of Kontsevich et al.,
at
least as exposed in [F], may be viewed as none other than the 'dual' of that of
[R]
for the case of rational curves, 'dual' meaning 'interchanging source and
target';
in particular, we sketch a proof from this viewpoint of the associativity of
quantum
multiplication.

This paper owes its existence to the unfailing encouragement of Bill Fulton,
who believed
all along in [R]; it is indeed a pleasure to thank him here.

\subheading{1.Old}
\vskip.10in
We find it technically convenient here to work with possibly reducible curves;
the modifications
or 'correction terms' needed to treat the irreducible case are a routine
matter.

Consider the locus $V_{d,\delta}$ of (not necessarily irreducible) curves of
degree
$d$ in $\Bbb P^2$ having $\delta$ ordinary nodes.  This is well known to be a
smooth
locally closed subvariety of pure codimension $\delta$ in $\Bbb
P^{\binom{d+2}{2}-1}$
and we are interested in its degree as such, which may be interpreted as the
number of
curves of $V_{d,\delta}$ passing through $\binom{d+2}{2} -\delta-1$ general
points in
$\Bbb P^2$, a number which we denote by $N_{d,\delta}$.  The idea is to get at
$N_{d,\delta}$
by a recursive procedure, based on specializing $\Bbb P^2$ to a surface (called
a `fan')
$$
S_0 = S_1 \cup S_2
$$
where $S_1 = Bl_0 (\Bbb P^2)$ (the `bottom' component), $S_2 = \Bbb P^2$ (the
`top'
component) and $E = S_1 \cap S_2$ (the 'axis') embedded in $S_i$ with
self-intersection $2i -3, i=1,2$.
Corresponding to this is a specialization
$$
V_{d,\delta} \to \sum m(\pi) V_{(d,e),(\delta_1, \delta_2),\pi},
\tag2
$$
where $V_{(d,e), (\delta_1, \delta_2),\pi}$ is a family of Cartier
divisors on $S_0$ whose general member $C_0$ may be described as follows:

$\bullet$\qquad $C_0 = C_1 \cup C_2$,

$\bullet$\qquad $C_1 \in |d H - eE|_{S_1}, C_2 \in |eE|_{S_2}$ nodal curves
with $\delta_1$
(resp. $\delta_2)$ nodes, smooth near $E$,

$\bullet$\qquad the divisor $D = C_1. E = C_2 . E$ has shape $\pi$, i.e. $\pi$
is a
patition having $\ell_i$ blocks of size $i$ (to be written as $\pi = [\ell_i])$
and $D =
\sum\limits_{i=1}^r \sum\limits_{j=1}^{\ell_i} i Q_{ij}, Q_{ij} \in E$
distinct.  Moreover
$m (\pi) = \prod\limits_{1}^r i^{\ell_i}$ and the sum is extended over all data
$((d,e), (\delta_1, \delta_2), \pi)$ satisfying
$$
\delta_1 + \delta_2 + \sum\limits_{i=1}^r (i-1) \ell_i = \delta
\tag3
$$
(i.e. each i-tacnode $iQ_{ij}$ `counts as $i-1$ nodes').

Now to apply the specialization (2) to the degree question, we specialize our
point set on
$\Bbb P^2$ to a collection of points on $S_0$, which a priori we may distribute
at will
among $S_1$ and $S_2$, with each distribution giving rise to some formula
which, however,
may or may not be usable.  For the purpose of the present discussion we will
make the
important simplifying assumption
$$
\delta < d ,
$$
and put $d+1$ points on $S_1$ and the remaining $\binom{d+1}{2} - 1 -\delta$ on
$S_2$.
It is then easy to see that the only limit components $V.$ that will contribute
to the
resulting formula will be ones with
$$
e = d-1.
$$
For those, we can write
$$
C_1 = C_{1,0} + \sum\limits_{i=1}^{\delta_1} R_i
\tag4
$$
with $C_{1,0}$ a smooth (rational) curve of `type' $(d-\delta_1, d-\delta_1 -
1)$
(i.e.  $C_{1,0} \in |(d-\delta_1) H - (d-\delta_1 - 1) E|)$ and $R_i$ distinct
rulings.

Now let us say that a partition $\pi' = [\ell_i'] \leq \pi = [\ell_i]$ if
$\ell_i' \leq
\ell_i  \forall i$, in which case we may define the complementary partition
$\pi - \pi' =
[\ell_i - \ell_i']$; also put $|\pi| = \sum i \ell_i, s(\pi) = \sum \ell_i, n
(\pi) =
\frac{s(\pi)!}{\ell_1! \cdots \ell_r !}$.  Counting the degree of a limit
component
$V_1 = \{C_1 \cup C_2\}$ in terms of those of $\{C_1\}$ and $\{C_2\}$ is
basically a matter
of decomposing the `diagonal' condition $C_1 . E = C_2 .E$ correspondingly to
the standard
Kunneth decomposition of the diagonal class on the product of $\Pi \Bbb
P^{\ell_i}$ with
itself; this leads to a sum of conditions corresponding to partitions $\pi'
\leq \pi$, each
amounting to fixing the location on $E$ of a portion $D'$ of $C_1 . E$
corresponding to
$\pi'$ and the complementary portion $D''$ of $C_2 . E$ corresponding to $\pi
-\pi'$.  The
resulting formula is as follows.
$$
\align
N_{d,\delta} &= \sum\limits_{|\pi| = d-1} m(\pi) \sum\Sb \pi' = [\ell_i']\\
\leq \pi =
[\ell_i]\endSb m(\pi - \pi') n(\pi - \pi')N_{d-1, \delta - s(\pi-\pi') + s(\pi)
-d+1, \pi -
\pi', \pi'}\\
&\times \sum\limits_{j=0}^{\ell_1'} {\binom{\ell_1'}{j}} {\binom{d+1}{s(\pi -
\pi') - j}} .
\tag5
\endalign
$$
Here $N_{e, \delta_2, \pi'',\pi'}$ denotes the degree of the locus of nodal
curves of
degree $e$  with $\delta_2$ nodes meeting a fixed line $E$ in a fixed divisor
of shape $\pi''$
plus a divisor of shape $\pi'$.  We have used the fact that $\delta_1  = s(\pi
-\pi')$, which
comes from the observation that the number of `axis' conditions on the bottom
curve $C_1$, i.e. $|\pi| - s(\pi - \pi') = d-1-s(\pi - \pi')$, plus the number
of
`interior' points imposed, i.e. $d+1$, must equal the dimension of the family
(4), i.e.
$2d-\delta_1$.  Also, the factor $m(\pi-\pi') n(\pi - \pi')$ is simply the
degree of the
`discriminant' variety of divisors of shape $\pi - \pi'$ on $E=\Bbb P^1$, while
the binomial
factors correspond to letting $j$ of the rulings go through some of the
multiplicity -1
part of $D'$ with the remaining $\delta_1 - j$ going through some of the $d+1$
interior
points.

Now of course in general the formula (5) is not by itself sufficient as one
needs a recursive
formula starting and ending with the $N_{e, \delta_2, \pi'',\pi'}$ or something
similar.
Such a formula is indeed given in [R], and it is not our purpose to reproduce
it here.  In
the examples worked out below the necessary further recursion is relatively
straightforward,
and will be indicated.

\subheading{Example 1}  $N_{d,2}$

There are seven relevant limit components and we proceed to list them and their
contributions.

A.  $V_{(d,d-1), (0,2),[d-1]}$; multiplicity $m=1$; contribution $N_{d-1,2}$

B.  $V_{(d,d-1),(1,1),[d-1]}$; $m =1$.  As $\delta_1 = 1$ we must take $\pi' =
[d-2], \pi -
\pi' = [1]$ so $j=0$ or 1 and the contribution is $(d+1+d-2).
N_{d-1,1,[1],[d-2]} = 3
(2d-1)(d-2)^2$.

C.  $V_{(d,d-1), (2,0), [d-1]}$; $m = 1$; $\pi' = [d-3]; j = 0,1,2$ ,
contribution $=
({\binom{d-3}{2}} + (d-3) (d+1) + {\binom{d+1}{2}} ). N_{d-1,0, [2],[d-3]} =
2d^2 - 5d+3$.

D.  $V_{(d,d-1),(0,1), [d-3,1]}$; $m = 2, \delta_1 = 0 \Rightarrow \pi' = \pi$,
so contribution
is $2N_{d-1,1,0,[d-3,1]}$. By an easier but simpler recursion (involving 1 node
and
1 tangency), the latter evaluates to $12(d-1)(d-2)(d-3)$.

E.  $V_{(d,d-1), (1,0),[d-3,1]}; m=2$.  $\pi' = [d-3]$ or $[d-4,1]$,
contribution
$= 8(d-1)(d-3)$.

F.  $V_{(d,d-1), (0,0),[d-4, 0,1]}, m = 3, \pi' = \pi$, contribution $9d-27$.

G.  $V_{(d,d-1), (0,0), [d-5,2]}, m = 4 , \pi' = \pi$, contribution 4.4.
${\binom{d-3}{2}}
= 9d^2 - 56 d + 96$.
\ls
Summing up, we get
$$
N_{d,2} - N_{d-1,2} = 18 d^3 - 81 d^2 + 84d + 12.
$$
Moreover it is easy to see that $N_{3,2} = {\binom{7}{2}} = 21$ so by
integrating we get
$$
N_{d,2} = \frac{9}{2} d^4 - 18 d^3 + 6 d^2 + \frac{81}{2} d - 33.
$$
This is a classical formula due to S. Roberts [Ro], which has been given modern
treatment
by I. Vainsencher [V].  Note that the curves are automatically irreducible if
$d \geq 4$.

\subheading{Example 2}  $N_{4,3}$

Here we have seven limit components.

A.  $V_{(4,3), (0,3), [3]}, m = 1$, contribution 15.

B.  $V_{(4,3), (1,2), [3]}, m =1$ contribution 21.7 = 147.

C.  $V_{(4,3), (2,1), [3]}, m = 1$, contribution $15 N_{3,1,[2],[1]} = 180$.

D.  $V_{(4,3), (3,0), [3]}, m = 1$, contribution ${\binom{5}{2}} = 10$.

E.  $V_{(4,3), (1,1), [1,1]}, m = 2, \pi' = [1]$ or $[0,1]$.  Contribution $2.2
N_{3,1,[0,1],[1]} + 2.5. N_{3,1,[1],[0,1]}$.  By a similar but simpler
recursion the
latter $N$'s evaluate respectively to 10, 16, so the total contribution is 200.

F.  $V_{(4,3), (0,2), [1,1]}, m = 2, \pi = \pi' = [1,1]$, contribution $2.15.2
= 60$.

G.  $V_{(4,3), (0,1), [0,0,1]}, m = 3, \pi = \pi' = [0,0,1]$, contribution
$3.N_{3,1,[0],[0,0,1]}$.  By a similar but simpler recursion, the latter $N$ is
21, so the contribution is 63.

Summing up, we get
$$
N_{4,3} = 675 = 5^2. 3^3.
$$

As the $\{{\text cubic + line}\}$ locus clearly has degree ${\binom{11}{2}} =
55$, we obtain 620 as the
number of irreducible rational quartics through 11 points.  (cf. [Z]).

\ls

\vskip .5in
\subheading{2.New}
\vskip.10in
The new approach works for maps from a fixed curve $C$, say to $\Bbb P ^2$. For
simplicity we
will assume $C=\Bbb P ^1$. Considering rational curves of degree $d$ in $\Bbb P
^2$
amounts to considering curves of bidegree $(1,d)$ in $\Bbb P ^1 \times \Bbb P
^2$, and
the old method to count them is by specialising the $\Bbb P ^2$ factor to a
fan;
the new approach on the other hand is to specialise the $\Bbb P ^1$ factor to a
'1-dimensional fan', i.e. to
$$
C_0 = C_1 \cup C_2, C_i = \Bbb P ^1, C_1 \cap C_2 = \{ x \}.
$$
Because $\Bbb P ^1$ is simpler than $\Bbb P ^2$ this approach works better in
this case; on the other hand it is apparently unknown how to make it
work when
the source curve is allowed to vary with moduli.

To be precise, fix a pair of points $y_1,y_2$ and a pair of lines $L_3,L_4$ in
$\P ^2$ and 4 points $x_1,...x_4 \in \P ^1 = C$ and consider curves of bidegree
$(1,d)$ in $C \times \P ^2$ containing $(x_1,y_1),(x_2,y_2)$ and meeting
$x_3 \times L_3, x_4 \times L_4 $, as well as a further collection of $3d-4$
'horizontal' lines
$C \times z_j$. We then specialise this to $C_0 \times \P ^2$ in two ways: (A)
$x_1...,x_4$ specialise to $x_{1,1},x_{2,1}\in C_1, x_{3,2},x_{4,2}\in C_2$;
(B) $x_1,x_3,x_2,x_4$ specialise to $x_{1,1},x_{3,1}\in C_1, x_{2,1}, x_{4,2}
\in C_2.$ In the (A) limit it is possible to have a component of bidegree
$(1,0)$ in $C_2 \times  (L_3 \cap L_4)$, while in the (B) limit all curves have
bidegrees $(1,d_1)\cup (1,d_2), d_1+d_2 = d, d_i > 0$. Thus letting $n_d$
denote the number of rational curves in $\P ^2$ through $3d-1$ points,
writing $(A) = (B)$ we get an equation of the form
$$
n_d + f(n_1,...,n_{d-1}) = g(n_1,...,n_{d-1})
$$
for suitable quadratic expressions $f,g$, which may be solved for
$n_d$.
\subheading{Example: $d=4$} $$ f = {\binom{8}{2}}.12.1.1.1.3 +
{\binom{8}{3}}.1.1.2.2.4 + 1.1.12.3.3.3 = 2228$$
with the summands corresponding to $d_1=3,2,1$ and, e.g. in the first product
the factors corresponding to: choosing 2 of the 8 points $z_j$ for
the image of $C_2$ to go through; the number of possible images of $C_1,
C_2, x_3,x_4,x;$
$$
g = 8.12.1.3.1.3 + {\binom{8}{4}}.1.1.2.2.4 + 8.1.12.1.3.3 = 2848
$$
$$
n_4 = 620.
$$
It is not difficult to parlay the above considerations into a general proof
of the main property of quantum cohomology, i.e. associativity. We now briefly
indicate how to do this, making free use of the notations of [F]
(and thus on the substantive level making use of the spaces
$\bar M_{0,n}(X,\beta )$, which might possibly be avoided as in the
example above). With $X$ being
as in [F] we may consider a space $\bar M_{0,0;n_1,n_2}(X,\beta _1,\beta _2)$
parametrising maps of $C_0$ into $X$ with $n_i$ points marked on $C_i, i=1,2$,
where
the image of $C_i$ has class $\beta _i,i=1,2$; this coincides with the
subvariety
of the product $\bar M_{0,n_1+1}(X,\beta_1)\times \bar M_{0,n_2+1}(X,\beta_2)$
given by the pullback of the diagonal $\Delta \subset X\times X$ via a suitable
projection.
Now fixing any decomposition of $n$ as $n_1+n_2$,
a divisor $D$ in the space $\bar M_{0,n+4}(X,\beta)$, defined by fixing
the position of the first four points, say, may be specialised to the union
(with
multiplicity 1) of the various $\bar M_{0,0;n_1+2,n_2+2}(X,\beta_1,\beta_2)$
with $\beta_1+\beta_2 = \beta$ (because fixing the position of three points
on each component of $C_0$ is vacuous); using the formula
$$
[\Delta] = \sum g^{e,f}T_e\otimes T_f,
$$
we conclude:
$$
I'_{\beta}(\gamma ^nT_iT_jT_kT_l) = \sum I_{\beta_1}(\gamma^{n_1}T_iT_jT_e)
I_{\beta_2}(\gamma^{n_2}T_kT_lT_f)g^{e,f}
$$
where the LHS is the appropriate Gromov-Witten integral over $D$ and
the sum is over all $\beta _1+\beta _2 = \beta$ and all $e,f$.
As the LHS is symmetric in $i,j,k,l$, this is certainly
stronger than the relation (2.9) on p.16 of [F], which is equivalent to
associativity.

{\bf References}
\roster
\item"{[D]}"  Dijkgraaf, R:  `Mirror symmetry and elliptic curves' (preprint).

\item"{[F]}"  Fulton, W:  `Enumerative geometry via quantum cohomology'
preprint.

\item"{[R]}"  Ran, Z.:  `Enumerative geometry of singular plane curves' Invent.
Math. {\bf 97}
(1989), 447-465.

\item"{[Ro]}"  Roberts, S.:  `Sur l'ordre des conditions $\ldots$' Crelle's J.
{\bf 67} (1867),
266-278.

\item"{[V]}"  Vainsencher, I.:  `Counting divisors with prescribed
singularities' Trans. AMS
{\bf 267}(1981), 399-422.

\item"{[Z]}"  Zeuthen, H.G., Pieri, M.:  `G\'eom\'etrie Enumerative.' Encyc.
Sci. Math. III 2,
260-331.  Leipzig:  Teubner 1915.
\endroster
\enddocument